\documentclass[11pt,a4paper,leqno]{article}
\usepackage{a4wide}
\usepackage{latexsym}
\usepackage{amsmath}
\usepackage{amssymb}
\usepackage{stackrel}  
\usepackage{color}
\usepackage{graphicx}
\usepackage[linesnumbered,ruled,vlined]{algorithm2e}

\usepackage{authblk}

\usepackage{hyperref}

\newtheorem{lemma}{Lemma}
\newtheorem{prop}{Proposition}
\newtheorem{theo}{Theorem}
\newtheorem{corol}{Corollary}
\newtheorem{defin}{Definition}

\newtheorem{remark}{Remark}

\newenvironment{proof}{\medskip\par\noindent{\bf Proof}}{\hfill $\Box$ \medskip\par}

\newcommand{\C}{\mathbb{C}}
\newcommand{\N}{\mathbb{N}}
\newcommand{\R}{\mathbb{R}}
\newcommand{\Z}{\mathbb{Z}}
\newcommand{\Oo}{\mathcal{O}}

\title{On integral representations of $q$-difference operators and their applications}

\begin{document}

\author[1]{Antonio C{\'a}ceres} 
\author[1]{Alberto Lastra} 
\author[2]{S{\l}awomir Michalik}
\author[2]{Maria Suwi{\'n}ska}
\affil[1]{Universidad de Alcal\'a, Dpto. F\'isica y Matem\'aticas, Alcal\'a de Henares, Madrid, Spain. {\tt antonio.caceresg@edu.uah.es, alberto.lastra@uah.es}
}
\affil[2]{Faculty of Mathematics and Natural Sciences, College of Science, Cardinal Stefan Wyszy{\'n}sky University, W\'oycickiego 1/3, 01-938 Warszawa, Poland\\
{\tt s.michalik@uksw.edu.pl, m.suwinska@op.pl}
}


\date{}

\maketitle
\thispagestyle{empty}
{ \small \begin{center}
{\bf Abstract}
\end{center}

Integral representations of two $q$-difference operators are provided in terms of special functions arising in the theory of asymptotic solutions to $q$-difference equations in the complex domain. Both representations are unified through the so-called $(p,q)$-differential operator, for which a kernel-like function is provided, generating the sequence of $(p,q)$-factorials.

\smallskip

\noindent Key words: integral representation, $q$-derivative, moment derivative, $(p,q)$-derivative. 

2020 MSC: 47G10, 44A20, 47A57, 44A35
}

\bigskip \bigskip

\footnotetext[0]{Acknowledgments: The second author is partially supported by the project PID2022-139631NB-I00 of Ministerio de Ciencia e Innovaci\'on, Spain.}

\section{Introduction}

The present work is devoted to the study of integral representations of various $q$-difference operators arising in the theory of functional equations in the complex domain, together with its application to the asymptotic representation of analytic functions in the Riemann surface of the logarithm, say $\mathcal{R}$, near the origin. 

The knowledge of asymptotic behavior of holomorphic functions defined in sectorial domains near the origin has been of great importance during the last decades. This phenomenon has been motivated by the possibility of giving an analytic meaning to the formal solutions to functional equations, initially to differential equations. We refer to~\cite{ba2} for a broader view on the topic. Previous asymptotic techniques have also been developed in the framework of $q-$difference equations, giving rise to a parallel theory since the end of the twentieth century, motivated by the $q-$Gevrey growth of the solutions to analytic linear $q-$difference equations determined in~\cite{bezivin}. In this direction, different $q$-analogs of the classical theory of Borel-Laplace summability procedure of formal power series arose.  

One of the choices for the construction of a Laplace operator leans on the $q$-difference equation satisfied by the Jacobi Theta function $\Theta_q$ (see Lemma~\ref{lema1}). Let us assume that $q\in\C$ with $|q|>1$. In~\cite{rz}, the authors' proposal is a discrete $q$-analog of Laplace transform, defined by
\begin{equation}\label{e0}
\mathcal{L}_{q;1}(g)(z)=\sum_{m\in\Z}\frac{g(q^m\lambda)}{\Theta_q\left(\frac{q^m\lambda}{z}\right)},
\end{equation}
for functions $g$, which is convergent near the origin and can be analytically prolonged on the discrete $q$-spiral $\lambda q^{\Z}=\{\lambda q^{m}:m\in\Z\}$ with $q$-exponential growth, i.e., they satisfy estimates of the form
\begin{equation}\label{e1}
|g(z)|\le|z|^{\mu}\exp\left(\frac{1}{2\log|q|}\log^2|z|\right),
\end{equation}
for some $\mu\in\R$. This technique has been successfully applied in the search of analytic and asymptotic solutions to concrete problems such as~\cite{d,lamasa12,taya}, among others. A continuous version of the previous $q$-Laplace transform was put forward in~\cite{z} and is defined by the following integral:
$$\mathcal{L}_{q;1}^d(g)(z)=\frac{1}{\log(q)}\int_{L_d}\frac{g(s)}{\Theta_q\left(\frac{s}{z}\right)}\frac{ds}{s},$$
for some $d\in\R$, where $L_d$ is a ray of complex numbers with argument $d$, departing from the origin to infinity, and where $g$ stands for a holomorphic function defined on some infinite sector with a bisecting direction $d$, subject to bounds as in (\ref{e1}) near infinity. This continuous version of $q$-Laplace transform has been applied in recent research such as~\cite{drlama,lama12,lama24,m23}.

On the other hand, a second branch of $q$-analogs have been considered in the literature, based on $q$-exponentials. In~\cite{vz}, the authors define two $q$-analogs of Laplace transform, a discrete and a continuous one, considering the $q$-exponential function
\begin{equation}\label{e2}
\hbox{exp}_q(z)=\sum_{n\ge0}\frac{1}{[n]_q^{!}}z^n
\end{equation}
as a kernel for the integral (continuous and discrete) operators defining different $q$-analogs of Laplace transforms:
$$\tilde{\mathcal{L}}^d_{q;1}(g)(z)=\frac{q-1}{\log(q)}\int_{L_d}\frac{g(t)}{\exp_q(qt/z)}dt,$$
and for $p=1/q$
$$\tilde{\mathcal{L}}_{q;1}(g)(z)=\frac{q}{1-1/q}\int_{\lambda p^{\Z}}\frac{g(t/(1-p))}{\exp_q(qt/((1-p)z))}d_pt,$$
where $\lambda p^{\mathbb{Z}}=\{\lambda p^{m}:m\in\Z\}$, and the previous integral stands for the discrete Jackson's $p$-integral (see~\cite{jackson}). In (\ref{e2}), the expression $[n]_q^{!}$ represents the $q$-factorial number (see Definition~\ref{defi137}). A second $q$-exponential function can also be defined, related to the $q$-Laplace transform in (\ref{e0}). See Definition~\ref{defi111}.

In~\cite{vz}, the previous definitions and those stated in terms of the Jacobi Theta function are compared.

Another noteworthy work is~\cite{tahara} by H. Tahara, where the adaptation of the previous operators considered by the author is based on the $q$-exponential function 
$$\hbox{Exp}_q(z)=\sum_{n\ge0}\frac{q^{\frac{n(n-1)}{2}}}{[n]_q^{!}}z^n,$$ 
which, contrary to the function from (\ref{e2}), is not an entire function, but is nevertheless convergent on some neighborhood of the origin. This new approach turns out to be more suitable for the study of certain families of $q$-difference equations, as can be seen in~\cite{tahara2}.

All the previous tools are to be considered for providing an integral representation to $q$-analogs of the derivative. Such $q$-analogs are interpreted as the realization of moment derivative when fixing two particular sequences, namely the sequence $(q^{\frac{n(n-1)}{2}})_{n\ge0}$ and $([n]_q^{!})_{n\ge0}$. A moment derivative is a formal operator associated with a fixed sequence of positive real numbers $m=(m_n)_{n\ge0}$, acting on the set of formal power series with complex coefficients, defined by
$$\partial_{m,z}\left(\sum_{n\ge0}\frac{u_n}{m_n} z^n\right)= \sum_{n\ge0}\frac{u_{n+1}}{m_n} z^n.$$
A theory of generalized summability and multisummability, associated with sequences $m$ has been recently developed in~\cite{sanz,lamasa,jisash}, searching for analytic and asymptotic solutions to moment-differential equations, i.e., functional equations under the action of moment derivatives. This theory applies to sequences $m$ which are strongly regular, leading to the existence of kernel functions for generalized summability, under broad assumptions. This summability theory has recently been applied to moment differential equations in research such as~\cite{lamisu,lamisu2,lamisu3,lamisu4,su}. The theory applies to sequences $m$ which are strongly regular. The sequences $([n]_q^{!})_{n\ge0}$ and $(q^{\frac{n(n-1)}{2}})_{n\ge0}$, with $q>1$, do not belong to this class of sequences. However, due to their importance in applications, different results have been achieved in these particular cases, as mentioned above. In the present work, we take into account the approach considered in the construction of convolution kernels for generalized summability in~\cite{jikalasa} and~\cite{ba2}, Section 5.8, to provide kernel functions for generalized summability associated with the sequence of $(p,q)$-factorials, $([n]_{p,q}^{!})_{n\ge0}$ which still remains outside the set of strongly regular sequences, where the theory of generalized summability remains valid. 
 
The $(p,q)$-derivative (see Definition~\ref{defi399}) can be also seen as a moment derivative associated with the sequence of $(p,q)$-factorials, $([n]_{p,q}^{!})_{n\ge0}$. Nowadays, this operator is of great importance as it has been applied in a wide variety of scientific fields, for example in quantum physics. We provide~\cite{kam,nun,prom} as recent examples of such a usage. The particularization to $(q,1)$-difference equations yields the study of $q$-difference equations, leading to more recent applications such as~\cite{prrasp,prsp}.

First two main results of the present work, Theorem~\ref{th:1} and Theorem~\ref{th:2}, provide integral representations of $q$-derivatives of functions holomorphic  near the origin by means of $q$-exponential functions and the Jacobi Theta function. These results are based on accurate estimations of the auxiliary special functions provided before. As a particularization of the previous result, one achieves two integral representations for the $(p,q)$-derivative (Corollary~\ref{coro1}).

In the last section, we also provide kernel functions associated with the sequence of $(p,q)$-factorials, taking into account that such sequence can be written as the product of two sequences. Both sequences do not satisfy the moderate growth condition (and therefore the theory of generalized summability does not apply). On the other hand, both sequences are sequences of moments, associated with known kernel functions. This fact allows us to construct a convolution kernel, associated with the sequence of $(p,q)$-factorials, in Theorem~\ref{th:3}. This is the key point for summability-like results, associated with such sequence, which is left to a separate research due to its intrinsic importance.

The paper is structured as follows. Section~\ref{sec2} recalls some of the main facts associated with different $q$-analogs of the exponential function, constructed by means of the moments associated with two kernel functions. The estimates associated with such functions will be crucial in Section~\ref{sec3} to provide integral representations of auxiliary functions in terms of $q$-exponentials and the Jacobi Theta function and to $q$-derivatives of  functions holomorphic near the origin in Theorem~\ref{th:1} and Theorem~\ref{th:2}. The last section, Section~\ref{sec6}, focuses on the $(p,q)$-derivative as a special situation, with great importance in applications. In the final result of the present work, Theorem~\ref{th:3}, we construct a convolution kernel associated with the sequence of $(p,q)$-factorials, which determines such derivative.

In the present study, we consider $q$ to be a real number with $q>1$. 

\section{Known facts on $q$-analogs}\label{sec2}

This section is devoted to the main facts about certain $q$-analogs in the complex domain. We recall the construction of two $q$-analogs of the exponential function which are established by means of the moments associated with certain positive measures with support in $[0,\infty)$. Known upper and lower bounds related to these elements will also be described.

We define the dilation operator $\sigma_q$ by 
$$\sigma_qf(z)=f(qz),$$
whenever the previous definition makes sense. Such operator is well defined on the set $\C[[z]]$, and also on functional spaces defined on sets $U$ such that $qU=\{qz:z\in U\}\subseteq U$. An example of such a set $U$ is any infinite sector in the complex domain with vertex at the origin.

Jackson's $q$-derivative is defined by
$$D_qf(z)=\frac{f(qz)-f(z)}{qz-z},\qquad z\neq0$$
on analogous sets.

We recall the definition of the Jacobi Theta function. 

\begin{defin}
\emph{The Jacobi Theta function} is defined by the formal power series
$$\Theta_{q}(z)=\sum_{n\in\Z}\frac{1}{q^{\frac{n(n-1)}{2}}}z^n,$$
determining a holomorphic function in $\C^{\star}=\C\setminus\{0\}$, and having an essential singularity at the origin.
\end{defin}

The proof of the following result can be found in~\cite{lama12}.

\begin{lemma}\label{lema1} The
Jacobi Theta function satisfies the following equality for every $m\in\Z$
\begin{equation}\label{e136}
\Theta_{q}(q^mz)=q^{\frac{m(m+1)}{2}}z^m\Theta_{q}(z),
\end{equation}
valid for all $z\in\C^{\star}$. In addition to this, for every $\delta>0$ there exists $\Delta=\Delta(q)>0$, which does not depend on $\delta$, such that 
$$\left|\Theta_{q}(z)\right|\ge \Delta\delta\exp\left(\frac{1}{2}\frac{\log^2|z|}{\log(q)}\right)|z|^{1/2},$$
valid for every element in
$$\C\setminus\left(\bigcup_{m\in\Z}\left\{z\in\C^{\star}: \left|1+zq^m\right|\le\delta\right\}\right).$$
\end{lemma}

%

We also recall several $q$-analogs of the exponential function appearing in the literature.

\begin{defin}\label{defi137}
Let $[0]_q=0$ and $[n]_{q}=\frac{q^{n}-1}{q-1}=\sum_{j=0}^{n-1}q^{j}$ for $n\ge 1$ be the $q$-numbers. We consider \emph{the $q$-factorials} defined by
$$
[n]_q^{!}=\begin{cases}
	1&\textrm{ for }n=0\\
	[1]_{q} [2]_{q}\cdots[n]_{q} &\textrm{ for }n\geq 1
\end{cases}
$$

The formal power series 
$$\exp_{q}(z)=\sum_{n\ge0}\frac{1}{[n]_q^!}z^n$$
defines an entire function.
\end{defin}

One can observe the confluence of the previous elements to the classical ones by taking $q\to1$: the $q$-numbers converge to the corresponding numbers. Consequently, also the $q$-factorials converge to the corresponding factorials, from which follows convergence of $\exp_{q}(z)$ to $\exp(z)$ for every $z\in\C$.

A second $q$-analog is linked to the theory of moment summability developed by W.~Balser (Section 5.5,~\cite{ba2}) and adapted to the generalized setting of strongly regular sequences by J.~Sanz in~\cite{sanz}. Although the sequence $(q^{\frac{n(n-1)}{2}})_{n\ge0}$ is not of moderate growth, one can adapt the~construction of the kernel functions for generalized summability in this context as follows.

\begin{defin}\label{defi111}
We define the formal power series
$$E_q(z)=\sum_{n\ge0}\frac{1}{q^{\frac{n(n-1)}{2}}}z^n,$$
which determines an entire function.
\end{defin}

We observe that (see (4.7) in~\cite{vz})
\begin{equation}\label{e169}
q^{\frac{n(n-1)}{2}}=\frac{q}{\log(q)}\int_{L_d}\frac{1}{\Theta_{q}(qt)}t^{n}dt,
\end{equation}
where the integral is performed along the infinite ray $L_d=[0,\infty)e^{id}$, for some fixed $d\in(-\pi,\pi)$. This entails that $(q^{\frac{n(n-1)}{2}})_{n\ge0}$ is the sequence of moments associated with the function $(\Theta_{q}(qz))^{-1}$, by choosing $d=0$. The application of such sequence in the theory of summability of formal solutions to functional equations is associated with a $q$-analog of Borel-Laplace procedure stated in~\cite{z} and successfully applied to solve concrete families of functional equations in~\cite{drlama,lama12,lamaq15,lasase,m23}, among others.
 
By virtue of Proposition 2.1~\cite{ramis}, one has the following result regarding the growth rate of $E_q$ at infinity.
\begin{lemma}\label{lema2}
There exist $K>0$ and $\alpha\in\R$ such that
$$|E_{q}(z)|\le K\exp\left(\frac{\log^2|z|}{2\log(q)}\right)|z|^{\alpha},$$
for every $z\in\C$
\end{lemma}

Some remarks are at hand at this point. First, we recall that the definition of $\exp_{q}(z)$ in Definition~\ref{defi137} is also determined by a generalized summability procedure. Indeed, regarding (4.5) in~\cite{vz}, one has
\begin{equation}\label{e176}
[n]_{q}^!=\frac{q-1}{\log(q)}\int_{L_d}t^n\exp_{1/q}(-qt)dt,
\end{equation}
for $d\in(-\pi,\pi)$. Therefore, $\exp_q(qz)$ can be seen as a kernel function for generalized summability associated with the sequence $([n]_q^!)_{n\ge0}$, which is not of moderate growth. We recall that $\exp_q(qt)\exp_{1/q}(-qt)\equiv 1$.

We also point out the relation between the sequences $(q^{\frac{n(n-1)}{2}})_{n\ge0}$ and $([n]_q^!)_{n\ge0}$, which satisfy
$$\lim_{n\to\infty}\frac{[n]_q^!}{q^{\frac{n(n-1)}{2}}}\left(\frac{q-1}{q}\right)^{n}=c(q),$$
for some constant $c=c(q)>0$ which only depends on $q$. Roughly speaking, both sequences are said to be equivalent in the sense that their growth at infinity is analogous up to a geometric term.

Moreover,  $[n]_q^!=[n]_{1/q}^!q^{\frac{n(n-1)}{2}}$ for $n\in\N_0$.
Hence, by Theorem~4~\cite{lami} the sequences $(q^{\frac{n(n-1)}{2}})_{n\ge0}$ and $([n]_q^!)_{n\ge0}$ have the same $q$-Gevrey asymptotic expansion in the sense that for every $s>0$ and $d\in\R$ the series $\sum_{n\ge0}\frac{a_n}{q^{\frac{n(n-1)}{2}}}t^n$ is convergent in a complex neighborhood of the origin 
 and its sum can be extended holomorphically on an infinite sector of bisecting direction $d$ with $q$-exponential growth of order $1/s$ on such sector if and only if the series $\sum_{n\ge0}\frac{a_n}{[n]_q^!}t^n$ has the same property.

A direct application of Proposition 2.1~\cite{tahara} or Lemma~2.6~\cite{zhang01}, leads to the following estimates for $\exp_q(z)$:

\begin{prop}\label{prop1}
\begin{itemize}
\item There exists $C_1>0$ such that 
$$|\exp_q(z)|\le C_1\exp\left(\frac{\log^2|z|}{2\log(q)}\right)|z|^{\frac{\log(q-1)}{\log(q)}-\frac{1}{2}},$$
for every $z\in\C$, $|z|\ge 1$. $\exp_q(z)$ is bounded from above in $D(0,1)$ as it is an entire function.
\item There exists $C_0>0$ such that $|\exp_q(z)|\ge C_0$ for all $z\in D(0,\frac{q^{1/2}}{q-1})$. Moreover, there exists $K_0>0$ such that for every sufficiently small $\epsilon>0$
$$ |\exp_q(z)|\ge\frac{\epsilon}{K_0}\exp\left(\frac{\log^2|z|}{2\log(q)}\right)|z|^{\frac{\log(q-1)}{\log(q)}-\frac{1}{2}},$$
for every $z\in\C\setminus\left(D(0,\frac{q^{1/2}}{q-1})\cup\bigcup_{m\ge0}\left\{z\in\C:\left|z+\frac{q^{m+1}}{q-1}\right|\le \frac{\epsilon q^{m+1}}{q-1}\right\}\right)$.
\end{itemize}
\end{prop}
\begin{proof}
The first part is a direct consequence of iv), Proposition 5.5~\cite{ramis} (see~\cite{wallisser}, p.~331). The second part is a direct consequence of Proposition 2.2~\cite{tahara}.
\end{proof}

\begin{remark} The zeros of $\exp_q(z)$ are easy to obtain and they are located at $-\frac{q^{m+1}}{q-1}$, for $m\in\N_0$.
\end{remark}

\section{Some integral representations associated with special functions}\label{sec3}

\begin{prop}\label{prop2}
For every $\omega,z\in\C^{\star}$, $\omega\neq z$, and $d\in\R$ such that $\arg(\omega/z)+d\neq\pm\pi$ and 
\begin{equation}\label{e230}
\left|\frac{\omega}{z}\right|> q^{\alpha-\frac{1}{2}},
\end{equation}
where $\alpha$ is determined in Lemma~\ref{lema2}, associated with the growth of the entire function $E_q$, one has
\begin{equation}\label{e202}
\frac{z}{\omega-z}=\frac{q}{\log(q)}\int_{L_d}\frac{E_q(\xi)}{\Theta_{q}\left(q\frac{\omega \xi}{z}\right)}d\xi,
\end{equation}
where $L_d=[0,\infty)e^{id}$.
\end{prop}
\begin{proof}
From a formal point of view, one has that 
\begin{align*}
\frac{q}{\log(q)}\int_{L_d}\frac{E_q(\xi)}{\Theta_q\left(q\frac{\omega \xi}{z}\right)}d\xi&=
\frac{q}{\log(q)}\int_{L_d}\sum_{n\ge0}\frac{1}{q^{\frac{n(n-1)}{2}}}\frac{\xi^n}{\Theta_q\left(q\frac{\omega \xi}{z}\right)}d\xi\\
&=\frac{q}{\log(q)}\sum_{n\ge0}\frac{1}{q^{\frac{n(n-1)}{2}}}\int_{L_d}\frac{\xi^n}{\Theta_q\left(q\frac{\omega \xi}{z}\right)}d\xi.
\end{align*}
From the change of variables $\omega\xi=zs$
one has that the previous expression equals
\begin{equation}\label{e233}
\frac{q}{\log(q)}\sum_{n\ge0}\left(\frac{z}{\omega}\right)^{n+1}\frac{1}{q^{\frac{n(n-1)}{2}}}\int_{L_{d'}}\frac{s^n}{\Theta_{q}(qs)}ds,
\end{equation}
with $d'=\arg(\omega/z)+d$. Regarding~(4.7) in Proposition 4.7~\cite{vz}, we receive
$$q^{\frac{n(n-1)}{2}}=\frac{q}{\log(q)}\int_{L_d}\frac{t^{n}}{\Theta_{q}(qt)}dt$$
for every $d'\in(-\pi,\pi)$ and all integers $n\ge0$. Moreover, the expression in~(\ref{e233}) turns into
$$\sum_{n\ge0}\left(\frac{z}{\omega}\right)^{n+1}=\frac{z}{\omega-z}.$$
This concludes the reasoning from a formal point of view.

From an analytic point of view, one can split the integral on the right-hand side of (\ref{e202}) into two parts: integrating in $L_{d,1}=[0,r_1]e^{id}$ for some fixed $r_1>0$, and in $L_{d,2}=[r_1,\infty)e^{id}$. After the change of variable $\omega\xi=zs$ one has that, in the case that $\arg(\omega/z)+d\neq\pm\pi$, then $s\mapsto E_{q}(zs/\omega)$ is bounded from above in $D\left(0,\frac{|\omega|}{|z|r_1}\right)$, and the integrability in $L_{d,1}$ is a consequence of (4.7) in Proposition 4.8,~\cite{vz}, for $n=0$. On the other hand, the integrability in $L_{d,2}$ is determined as follows. 

From Lemma~\ref{lema1} and Lemma~\ref{lema2} we conclude that there exist $K>0$ and $\alpha\in\R$ and for any fixed $\delta>0$ there exists $\Delta>0$, such that
\begin{align*}
\left|\frac{E_q(\xi)}{\Theta_{q}\left(q\frac{\omega\xi}{z}\right)}\right|&\le 
\frac{K}{\Delta\delta}\exp\left(\frac{\log^2|\xi|-\log^2\left|q\frac{\omega\xi}{z}\right|}{2\log(q)}\right)|\xi|^{\alpha-\frac{1}{2}}\left|\frac{z}{q\omega}\right|^{\frac{1}{2}}\\
&=\frac{K}{\Delta\delta}\exp\left(\frac{-1}{2\log(q)}(\log^2\left|q\frac{\omega}{z}\right|+2\log|\xi|\log\left|q\frac{\omega}{z}\right|)\right)|\xi|^{\alpha-\frac{1}{2}}\left|\frac{z}{q\omega}\right|^{\frac{1}{2}}
\end{align*}
for every $z,\omega\in\C$ such that 
$$\left|1+\frac{\omega}{z}q^{m}re^{id}\right|>\delta,\quad\hbox{for every }m\in\Z,\ r>r_2.$$
We observe that the set of complex numbers satisfying the previous condition is not empty, provided that $\delta>0$ is sufficiently small and $d+\arg(\omega/z)\neq\pm\pi$.

The previous estimates and the choice (\ref{e230}) yield
\begin{equation}
\label{eq:E_q}
\left|\frac{E_q(\xi)}{\Theta_{q}\left(q\frac{\omega\xi}{z}\right)}\right|\le \frac{K}{\Delta\delta}\exp\left(\frac{-1}{2\log(q)}\log^2\left|q\frac{\omega}{z}\right|\right)\left|\frac{z}{q\omega}\right|^{\frac{1}{2}}|\xi|^{\alpha-\frac{1}{2}-1-\frac{1}{\log(q)}\log\left|\frac{\omega}{z}\right|}.
\end{equation}
The condition~(\ref{e230}) provides convergence of the integral in $L_{d,2}$.
\end{proof}

\begin{prop}\label{prop3}
For every $\omega,z\in\C^{\star}$, and $d\in\R$ such that $\arg(\omega/z)+d\neq\pm\pi$ and $\left|\frac{\omega}{z}\right|>1,$
one has
\begin{equation}\label{e203}
\frac{z}{\omega-z}=\frac{q-1}{\log(q)}\int_{L_d}\frac{\exp_q(\xi)}{\exp_q\left(q\frac{\omega \xi}{z}\right)}d\xi,
\end{equation}
where $L_d=[0,\infty)e^{id}$.
\end{prop}
\begin{proof}
From a formal point of view, one has that 
\begin{align*}
\frac{q-1}{\log(q)}\int_{L_d}\frac{\exp_q(\xi)}{\exp_q\left(q\frac{\omega \xi}{z}\right)}d\xi&=
\frac{q-1}{\log(q)}\int_{L_d}\sum_{n\ge0}\frac{1}{[n]_q^!}\frac{\xi^n}{\exp_q\left(q\frac{\omega \xi}{z}\right)}d\xi\\
&=\frac{q-1}{\log(q)}\sum_{n\ge0}\frac{1}{[n]_q^!}\int_{L_d}\frac{\xi^n}{\exp_q\left(q\frac{\omega \xi}{z}\right)}d\xi.
\end{align*}
From the change of variables $\omega\xi=zs$ it follows that the previous expression equals
\begin{equation}\label{e234}
\frac{q-1}{\log(q)}\sum_{n\ge0}\left(\frac{z}{\omega}\right)^{n+1}\frac{1}{[n]_{q}^!}\int_{L_{d'}}\frac{s^n}{\exp_q(qs)}ds,
\end{equation}
with $d'=\arg(\omega/z)+d$. Regarding~(4.5) in Proposition 4.7~\cite{vz} we get that for every $d'\in(-\pi,\pi)$ and all integer $n\ge0$ (see~\ref{e176}), the expression in (\ref{e234}) turns into
$$\sum_{n\ge0}\left(\frac{z}{\omega}\right)^{n+1}=\frac{z}{\omega-z}.$$
This concludes the reasoning from a formal point of view.

For the analytic proof of the result, we make use of Proposition~\ref{prop1}. Let $\epsilon>0$ be sufficiently small, and let us choose $z,\omega\in\C^{\star}$ with $z\neq\omega$ and such that
$$\left|1+r\frac{e^{id}\omega}{z}\right|\ge \epsilon,\qquad \hbox{for all }r\ge0.$$
This condition holds for small enough $\epsilon$ if $\arg(\omega/z)+d\neq\pm\pi$. Hence, we write 
$$\int_{L_d}\frac{\exp_q(\xi)}{\exp_q\left(q\frac{\omega \xi}{z}\right)}d\xi=\int_{L_{d,1}}\frac{\exp_q(\xi)}{\exp_q\left(q\frac{\omega \xi}{z}\right)}d\xi+\int_{L_{d,2}}\frac{\exp_q(\xi)}{\exp_q\left(q\frac{\omega \xi}{z}\right)}d\xi,$$
where $L_{d,1}=\left[0,\left|\frac{z}{w}\right|\frac{1}{q^{1/2}(q-1)}\right]e^{id}$, and $L_{d,2}=\left[\left.\left|\frac{z}{w}\right|\frac{1}{q^{1/2}(q-1)},\infty\right.\right)e^{id}$.
Now let us notice that, from the first part of Proposition~\ref{prop1}, we have
\begin{align*}
\left|\int_{L_{d,1}}\frac{\exp_q(\xi)}{\exp_q\left(q\frac{\omega \xi}{z}\right)}d\xi\right|&\le \frac{1}{C_0}\max_{0\le r\le \frac{1}{q^{1/2}(q-1)}\left|\frac{z}{\omega}\right|}|\exp_{q}(re^{id})|\left|\frac{z}{\omega}\right|\frac{1}{q^{1/2}(q-1)}\\
&=\frac{1}{C_0}\exp_q\left(\frac{1}{q^{1/2}(q-1)}\left|\frac{z}{\omega}\right|\right)\left|\frac{z}{\omega}\right|\frac{1}{q^{1/2}(q-1)}.
\end{align*}
On the other hand, let $\epsilon$ be a small positive number. From the second part of Proposition~\ref{prop1}, it holds that
\begin{align*}
\left|\int_{L_{d,2}}\frac{\exp_q(\xi)}{\exp_q\left(q\frac{\omega \xi}{z}\right)}d\xi\right|&
\le \frac{C_1K_0}{\epsilon}q^{\frac{1}{2}-\frac{\log(q-1)}{\log(q)}}\left|\frac{\omega}{z}\right|^{\frac{1}{2}-\frac{\log(q-1)}{\log(q)}}\\
&\times\int_{\left|\frac{z}{\omega}\right|\frac{1}{q^{1/2}(q-1)}}^{\infty}r^{\frac{\log(q-1)}{\log(q)}-\frac{1}{2}-\frac{\log(q-1)}{\log(q)}+\frac{1}{2}}\exp\left(\frac{\log^2(r)}{2\log(q)}-\frac{\log^2\left|\frac{qr\omega}{z}\right|}{2\log(q)}\right)dr
\end{align*}
whenever 
$$\left|\frac{re^{id}\omega}{z}+\frac{q^m}{q-1}\right|\ge \frac{\epsilon q^{m}}{q-1},$$
for every integer $m>0$ and all $r>\frac{1}{q^{1/2}(q-1)}\left|\frac{z}{\omega}\right|$. In particular, this condition holds for all $z,\omega$ such that 
$$\left|1+r\frac{e^{id}\omega}{z}\right|\ge \epsilon,\qquad \hbox{for all }r\ge0.$$
The previous estimates equal
\begin{multline*}
\frac{C_1K_0}{\epsilon}q^{\frac{1}{2}-\frac{\log(q-1)}{\log(q)}}\left|\frac{\omega}{z}\right|^{\frac{1}{2}-\frac{\log(q-1)}{\log(q)}}\exp\left(-\frac{\log^2(q\left|\frac{\omega}{z}\right|)}{2\log(q)}\right)\hfill
\\
\hfill\times\int_{\left|\frac{z}{\omega}\right|\frac{1}{q^{1/2}(q-1)}}^{\infty}\exp\left(-\frac{2\log(r)\log(q\left|\frac{\omega}{z}\right|)}{2\log(q)}\right)dr\\
=\frac{C_1K_0}{\epsilon}q^{\frac{1}{2}-\frac{\log(q-1)}{\log(q)}}\left|\frac{\omega}{z}\right|^{\frac{1}{2}-\frac{\log(q-1)}{\log(q)}}\exp\left(-\frac{\log^2(q\left|\frac{\omega}{z}\right|)}{2\log(q)}\right)\int_{\left|\frac{z}{\omega}\right|\frac{1}{q^{1/2}(q-1)}}^{\infty}r^{-\frac{\log\left(q\left|\frac{\omega}{z}\right|\right)}{\log(q)}}dr
\end{multline*}

We observe that the previous integral is convergent for 
$$-\frac{\log\left(q\left|\frac{\omega}{z}\right|\right)}{\log(q)}<-1,\quad\textrm{or equivalently}\quad\left|\frac{\omega}{z}\right|>1.$$
\end{proof}

\section{Integral representations of $q$-derivatives}\label{sec4}
\begin{defin}
Let $\mathbb{E}$ be a complex Banach space. For a given sequence of positive real numbers $m=(m(n))_{n\ge0}$ with $m(0)=1$, the formal operator
$\partial_{m,z}\colon\mathbb{E}[[z]] \to \mathbb{E}[[z]]$ satisfying
\begin{equation*}
\partial_{m,z}\left(\sum_{n\ge0}a_nz^n\right)=\sum_{n\ge0}\frac{m(n+1)}{m(n)}a_{n+1}z^n
\end{equation*}
or equivalently
\begin{equation*}
\partial_{m,z}\left(\sum_{n\ge0}\frac{u_n}{m(n)}z^n\right)=\sum_{n\ge0}\frac{u_{n+1}}{m(n)}z^n
\end{equation*}
is called the \emph{moment derivative}.
\end{defin}

Applying the operator $\partial_{m,z}$ to the Maclaurin series expansions of holomorphic functions,
we conclude that the moment derivative $\partial_{m,z}$
is well-defined on the space $\Oo(D)$ of holomorphic functions in a complex neighborhood of the origin.

\begin{remark}\label{rem:sequences}
Let us assume that $f(z)=\sum_{n\geq 0}a_nz^n\in\Oo(D)$. Note that for some special sequences of positive numbers $m=(m(n))_{n\ge 0}$ the moment derivatives $\partial_{m,z}$ become the well-known operators acting on $f(z)$. 
 \begin{enumerate}
  \item If $m=(n!)_{n\ge0}$ then $\partial_{m,z}f(z)=\sum_{n\ge0}(n+1)a_{n+1}z^n=f'(z)$ is the ordinary derivative.
  \item If $m_1=\left(q^{\frac{n(n-1)}{2}}\right)_{n\geq 0}$ then
  $\partial_{m_1,z}f(z)=\sum_{n\ge0}q^na_{n+1}z^n=\frac{f(qz)-f(0)}{qz}=:\tilde{D}_qf(z)$.
  \item If $m_2=([n]_q^!)_{n\ge0}$ then $\partial_{m_2,z}f(z)=\sum_{n\ge0}[n+1]_qa_{n+1}z^n=\frac{f(qz)-f(z)}{qz-z}=D_qf(z)$ is Jackson's $q$-derivative.
 \end{enumerate}
\end{remark}

From this point forward, we will mainly focus on the last two examples $\partial_{m_1,z}$ and $\partial_{m_2,z}$ of the moment derivatives.

Similarly to Proposition~3~\cite{mi}, we may construct the integral representations of $\partial_{m_1,z}$ and $\partial_{m_2,z}$ moment derivatives of the holomorphic function $f\in\Oo(D)$. To this end we use the Cauchy integral formula with the Cauchy kernel given by Propositions~\ref{prop2} and~\ref{prop3}, respectively.
Namely for $\tilde{D}_q$ derivative we get

\begin{theo}\label{th:1}
 Let $f\in\Oo(D(0,r))$ for some $r>0$. Then for every $\varepsilon\in(0,r)$, $k\in\N_0$ and $|z|<q^{-k}\varepsilon\min\{1,q^{\frac{1}{2}-\alpha}\}$, where $\alpha$ is determined in Lemma \ref{lema2}, we get
 \begin{equation}
  \tilde{D}^k_qf(z)=D_{0,q}^kf(z)=\partial_{m_1,z}^kf(z)=\frac{q}{2\log(q)\pi i}\oint_{|\omega|=\varepsilon}f(\omega)\int_{L_{\theta}}\zeta^k\frac{E_q(z\zeta)}{\Theta_q(q\omega\zeta)}\,d\zeta\,d\omega,
 \end{equation}
where $L_{\theta}=[0,\infty)e^{i\theta}$ and $\arg(\omega)+\theta\neq \pm\pi$.
\end{theo}
\begin{proof}
Applying the Cauchy Integral Formula and Proposition \ref{prop2} to $f\in\Oo(D(0,r))$ with $r>\varepsilon>0$, we observe that for every $|z|<\varepsilon\min\{1,q^{\frac{1}{2}-\alpha}\}$
\begin{equation*}
 f(z)=\frac{1}{2\pi i}\oint_{|\omega|=\varepsilon}f(\omega)\frac{z}{\omega-z}\frac{d\omega}{z}=\frac{q}{2\log(q)\pi i}\oint_{|\omega|=\varepsilon}f(\omega)\int_{L_d}\frac{E_q(\xi)}{\Theta_q\left(q\frac{\omega\xi}{z}\right)}\frac{d\xi}{z}\,d\omega
\end{equation*}
for $\arg(\omega)-\arg(z)+d\neq\pm\pi$.
After the change of variables $\zeta=\frac{\xi}{z}$ we conclude that
\begin{equation*}
f(z)=\frac{q}{2\log(q)\pi i}\oint_{|\omega|=\varepsilon}f(\omega)\int_{L_{\theta}}\frac{E_q(z\zeta)}{\Theta_q(q\omega\zeta)}\,d\zeta\,d\omega,
\end{equation*}
where $\theta=d-\arg(z)$. Hence $\arg(\omega)+\theta\neq\pm\pi$.

Calculating $\partial_{m_1,z}$ moment derivatives of the function $z\mapsto E_q(z\zeta)$ for fixed $\zeta\in\C$ we obtain
\begin{equation*}
 \partial_{m_1,z}^k E_q(z\zeta)=\partial_{m_1,z}^k\left(\sum_{n\ge 0}\frac{\zeta^n}{q^{\frac{n(n-1)}{2}}}z^n\right)=\sum_{n\ge 0}\frac{\zeta^{k+n}}{q^{\frac{n(n-1)}{2}}}z^n=\zeta^k\left(\sum_{n\ge 0}\frac{\zeta^n}{q^{\frac{n(n-1)}{2}}}z^n\right)=\zeta^k E_q(z\zeta).
\end{equation*}
Hence for every $z\in D(0,r)$ satisfying $q^k|z|<|\omega|=\varepsilon$ we get
\begin{multline*}
\partial_{m_1,z}^kf(z)=\frac{q}{2\log(q)\pi i}\oint_{|\omega|=\varepsilon}f(\omega)\int_{L_{\theta}}\frac{\partial_{m_1,z}^kE_q(z\zeta)}{\Theta_q(q\omega\zeta)}\,d\zeta\,d\omega\\
=\frac{q}{2\log(q)\pi i}\oint_{|\omega|
=\varepsilon}f(\omega)\int_{L_{\theta}}\zeta^k\frac{E_q(z\zeta)}{\Theta_q(q\omega\zeta)}\,d\zeta\,d\omega.
\end{multline*}
Moreover, by (\ref{eq:E_q}) with $\xi=z\zeta$ we get
\begin{equation*}
 \left|\zeta^k\frac{E_q(z\zeta)}{\Theta_q(q\omega\zeta)}\right| \leq
 \frac{K}{\Delta\delta}\exp\left(\frac{-1}{2\log(q)}\log^2\left|q\frac{\omega}{z}\right|\right)\left|\frac{z}{q\omega}\right|^{\frac{1}{2}}|z|^{\alpha-\frac{1}{2}-1-\frac{1}{\log(q)}\log\left|\frac{\omega}{z}\right|}|\zeta|^{k+\alpha-\frac{1}{2}-1-\frac{1}{\log(q)}\log\left|\frac{\omega}{z}\right|}.
\end{equation*}
Therefore the integral in $L_{\theta}$ is convergent for $\left|\frac{z}{\omega}\right|<q^{\frac{1}{2}-\alpha-k}$.
\end{proof}
\bigskip

We receive an analogous theorem for $D_q$ derivative:
\begin{theo}\label{th:2}
 Let $f\in\Oo(D(0,r))$ for some $r>0$. Then for every $\varepsilon\in(0,r)$, $k\in\N_0$ and $|z|<q^{-k}\varepsilon$,  we get
 \begin{equation}
  D^k_qf(z)=D^k_{1,q}f(z)=\partial_{m_2,z}^kf(z)=\frac{q-1}{2\log(q)\pi i}\oint_{|\omega|=\varepsilon}f(\omega)\int_{L_{\theta}}\zeta^k\frac{\exp_q(z\zeta)}{\exp_q(q\omega\zeta)}\,d\zeta\,d\omega,
 \end{equation}
where $L_{\theta}=[0,\infty)e^{i\theta}$ and $\arg(\omega)+\theta\neq\pm\pi$.
\end{theo}
\begin{proof}
Similarly to the proof of Theorem~\ref{th:1} we use the Cauchy Integral Formula and Proposition~\ref{prop3} to $f\in\Oo(D(0,r))$ with $r>\varepsilon>0$. Then for every $|z|<\varepsilon$ we get
\begin{equation*}
 f(z)=\frac{1}{2\pi i}\oint_{|\omega|=\varepsilon}f(\omega)\frac{z}{\omega-z}\frac{d\omega}{z}=\frac{q-1}{2\log(q)\pi i}\oint_{|\omega|=\varepsilon}f(\omega)\int_{L_d}\frac{\exp_q(\xi)}{\exp_q\left(q\frac{\omega\xi}{z}\right)}\frac{d\xi}{z}\,d\omega
\end{equation*}
for $\arg(\omega)-\arg(z)+d\neq\pm\pi$.
After performing the change of variables $\zeta=\frac{\xi}{z}$ we conclude that
\begin{equation*}
f(z)=\frac{q-1}{2\log(q)\pi i}\oint_{|\omega|=\varepsilon}f(\omega)\int_{L_{\theta}}\frac{\exp_q(z\zeta)}{\exp_q(q\omega\zeta)}\,d\zeta\,d\omega,
\end{equation*}
where $\theta=d-\arg(z)$ with $\arg(\omega)+\theta\neq\pm\pi$.

Calculating $\partial_{m_2,z}$ moment derivatives of the function $z\mapsto \exp_q(z\zeta)$ for fixed $\zeta\in\C$ we obtain
\begin{equation*}
 \partial_{m_2,z}^k \exp_q(z\zeta)=\partial_{m_2,z}^k\left(\sum_{n\ge 0}\frac{\zeta^n}{[n]_q^!}z^n\right)=\sum_{n\ge 0}\frac{\zeta^{k+n}}{[n]_q^!}z^n=\zeta^k\left(\sum_{n\ge 0}\frac{\zeta^n}{[n]_q^!}z^n\right)=\zeta^k \exp_q(z\zeta).
\end{equation*}
Hence, for every $z\in D(0,r)$ satisfying $q^k|z|<|\omega|=\varepsilon$ we get
\begin{multline*}
\partial_{m_2,z}^kf(z)=\frac{q-1}{2\log(q)\pi i}\oint_{|\omega|=\varepsilon}f(\omega)\int_{L_{\theta}}\frac{\partial_{m_2,z}^k\exp_q(z\zeta)}{\exp_q(q\omega\zeta)}\,d\zeta\,d\omega\\
=\frac{q-1}{2\log(q)\pi i}\oint_{|\omega|
=\varepsilon}f(\omega)\int_{L_{\theta}}\zeta^k\frac{\exp_q(z\zeta)}{\exp_q(q\omega\zeta)}\,d\zeta\,d\omega.
\end{multline*}
Repeating the proof of Proposition \ref{prop3} with $\xi=z\zeta$ we conclude that the integral
\begin{equation*}
 \int_{L_{\theta}}\zeta^k\frac{\exp_q(z\zeta)}{\exp_q(q\omega\zeta)}\,d\zeta
\end{equation*}
is convergent for
$$
k-\frac{\log\left(q\left|\frac{\omega}{z}\right|\right)}{\log(q)}<-1,\quad\textrm{or equivalently}\quad\left|\frac{\omega}{z}\right|>q^k.$$
\end{proof}

\section{A unified approach: the $(p,q)$-derivative}\label{sec6}

In this section, we provide a unified approach to the $q$-derivatives appearing in the previous section using post-quantum calculus (see \cite{gura}) and $(p,q)$-derivatives. Such derivatives turn out to be the moment derivatives associated with the sequence of $(p,q)$-factorials. We provide kernel functions generating such sequence, as the convolution of two known kernel functions.

\begin{defin}\label{defi399}
 Let $p\neq q$. The \emph{$(p,q)$-derivative} of the function $f$ is defined as
\begin{equation*}
 D_{p,q}f(z)=\frac{f(pz)-f(qz)}{pz-qz},\qquad z\neq 0.
\end{equation*}
\end{defin}

Putting $p=0$ and $p=1$ respectively in the above definition of the $(p,q)$-derivative $D_{p,q}$ we get the moment derivatives $\partial_{m_1,z}$ and $\partial_{m_2,z}$, respectively (see Section~\ref{sec4}).

We may extend the integral representations of $D^k_{0,q}f(z)$ and $D^k_{1,q}f(z)$ given in Theorems \ref{th:1} and \ref{th:2} to any $(p,q)$-derivative of order $k$ of the function $f$, where $0\le p < q$.
To this end it is sufficient to observe that for any $0<p<q$ and any $k\in\N_0$ holds
\begin{equation*}
 D^k_{p,q}f(z)=\frac{1}{p^k}D^k_{1,q/p}f(\tilde{z})\bigg|_{\tilde{z}=p^kz}.
\end{equation*}

\begin{corol}\label{coro1}
 Let $0\leq p <q$, $f\in\Oo(D(0,r))$ for some $r>0$, $\varepsilon\in(0,r)$ and $k\in\N_0$. Then $D^k_{p,q}f(z)$ has the following integral representations:
 \begin{enumerate}
  \item If $p=0$ then for every $|z|<q^{-k}\varepsilon\min\{1,q^{\frac{1}{2}-\alpha}\}$ we have
 \begin{equation*}
  D_{p,q}^kf(z)=\frac{q}{2\log(q)\pi i}\oint_{|\omega|=\varepsilon}f(\omega)\int_{L_{\theta}}\zeta^k\frac{E_q(z\zeta)}{\Theta_q(q\omega\zeta)}\,d\zeta\,d\omega.
 \end{equation*}
 \item If $p>0$ then for every $|z|<q^{-k}\varepsilon$ we get
 \begin{equation*}
  D^k_{p,q}f(z)=\frac{q/p-1}{2p^k\log(q/p)\pi i}\oint_{|\omega|=\varepsilon}f(\omega)\int_{L_{\theta}}\zeta^k\frac{\exp_{q/p}(p^kz\zeta)}{\exp_{q/p}(\frac{q}{p}\omega\zeta)}\,d\zeta\,d\omega.
  \end{equation*}
 \end{enumerate}
In both cases $L_{\theta}=[0,\infty)e^{i\theta}$ and $\arg(\omega)+\theta\neq \pm\pi$.
\end{corol}

\begin{remark}
 Since $D_{p,q}f(z)=D_{q,p}f(z)$, by the above corollary we obtain the integral representation of $D_{p,q}^kf(z)$ for every $p,q\in[0,\infty)$, $p\neq q$.
\end{remark}

It is straightforward to check that, given $p\neq q$, then the $(p,q)$-derivative of a formal power series corresponds to the moment derivative associated with the sequence of $(p,q)$-factorials, $m_{p,q}=([n]_{p,q}^{!})_{n\ge0}$. This sequence is defined as follows
$$
[n]_{p,q}^{!}=\begin{cases}
	1&\textrm{ for }n=0\\
	\displaystyle\prod_{j=1}^{n}\frac{p^j-q^j}{p-q}&\textrm{ for }n\ge 1
\end{cases}.
$$
From now on, we assume that $p$ and $q$ are positive real numbers.

\begin{lemma}\label{lema482}
Let $p,q>0$. Then, for every $n\ge0$ one has  
$$[n]_{p,q}^{!}=[n]_{p/q}^{!}q^{\frac{n(n-1)}{2}}\quad\textrm{for } n\ge0.$$
\end{lemma}
\begin{proof}
It is straightforward to check the identity for $n=0$. Observe that for all $1\le j\le n$ one has that 
$$\frac{p^j-q^j}{p-q}=\frac{(p/q)^j-1}{p/q-1}q^{j-1}.$$
The conclusion is a direct consequence of the definition of the $p/q$-factorial.
\end{proof}

As a consequence of Lemma~\ref{lema482}, the sequence of $(p,q)$-factorials can be written as
$$m_{p,q}=m_{p/q}m_1,$$
where $m_{p/q}$ is the sequence $m_2$ of Section~\ref{sec4}, associated with Jackson's $p/q$-derivative, by substituting $q$ by $p/q$, and $m_1$ is the sequence appearing in point~2 of Remark~\ref{rem:sequences}. 

Seeing as $m_{p,q}$ is a sequence of moments which can be written as the product of two sequences of moments, it is natural to consider the construction of a kernel function associated with it, in terms of each factor. This was previously done in Section 5.8~\cite{ba2}, by defining the convolution of kernel functions which determine a new kernel whose sequence of moments is the product of the initial sequences. However, this approach is considered whenever both sequences are equivalent to a Gevrey sequence. In~\cite{jikalasa}, the construction of convolution kernels is given for the larger class of moment sequences associated with a pair of kernel functions for generalized summability which admit a nonzero proximate order. It turns out that the sequence $m_{p,q}$ is not of moderate growth whenever $p$ or $q$ are larger than $1$, and the previous results do not apply in this situation, which is very significant for various applications. Indeed, if one of the numbers, say $q$, is smaller than 1, then the growth rate of the sequence $m_{p,q}$ is analogous to that of  the sequence $([n]_{p}^{!})_{n\ge0}$. Therefore, one can restrict the study to $p>q>1$.
 
In view of~(\ref{e176}) and~(\ref{e169}), the sequences $m_{p/q}$ and $m_1$ are sequences of moments associated with the kernel functions 
$$e_1(z)=\frac{p/q-1}{\log(p)-\log(q)}\frac{1}{\exp_{p/q}(\frac{pz}{q})},$$
and
$$e_2(z)=\frac{q}{\log(q)}\frac{1}{\Theta_{q}(qz)},$$
respectively. 

\begin{remark} It is worth mentioning that, in order that the kernel functions and the associated sequences of moments follow the construction by W.~Balser in~\cite{ba2}, a shift between the sequences of moments and the moments associated with the kernel functions should appear. Therefore, one has
$$[n]_{p/q}^{!}=\int_{0}^{\infty}t^{n-1}\tilde{e}_1(t)dt,\qquad n\ge1,$$
where $\tilde{e}_1(z)=ze_1(z)$, and 
$$q^{\frac{n(n-1)}{2}}=\int_{0}^{\infty}t^{n-1}\tilde{e}_2(t)dt,\qquad n\ge1,$$
where $\tilde{e}_2(z)=ze_2(z)$.
\end{remark}

In the next result, we describe the convolution of the previous kernels, which determine a kernel function associated with the product sequence of $(p,q)$-factorials, $m_{p,q}$.

The procedure is heavily based on the one described in Section~5.8~\cite{ba2} as well as Section~4.3~\cite{jikalasa}, so we omit details that remain unchanged.

The following classical definitions and results can be adapted to our settings. 

\begin{defin}
Let $\mathbb{M}=(M_n)_{n\ge0}$ be a sequence of positive real numbers. We say that $\mathbb{M}$ is a weight sequence if it is logarithmically convex, i.e., $M_n^2\le M_{n-1}M_{n+1}$ for every $n\ge1$, and $\lim_{n\to\infty}M_n^{1/n}=+\infty$.
In association with a weight sequence $\mathbb{M}$, one can define the function $\omega_{\mathbb{M}}$ given by $\omega_{\mathbb{M}}(0)=0$, and
$$\omega_{\mathbb{M}}(t)=\sup_{n\ge0}\log\left(\frac{t^n}{M_n}\right),\qquad t>0.$$
\end{defin}

It is straightforward to check that the sequences $([n]_{p/q}^{!})_{n\ge0}$ and $(q^{\frac{n(n-1)}{2}})_{n\ge0}$ are both well defined weight sequences. 

\begin{lemma}[Lemma 4.13,~\cite{jikalasa}]\label{lema4}
Given two weight sequences, $\mathbb{M}_j$ for $j=1,2$, then for all $r,s>0$ one has that
$$e^{\omega_{\mathbb{M}_1\cdot\mathbb{M}_2}(r)}\le e^{\omega_{\mathbb{M}_1}(s)}e^{\omega_{\mathbb{M}_2}(r/s)},$$
where the $n$-th element of the sequence $\mathbb{M}_1\cdot\mathbb{M}_2$ is given by the product of the $n$-th elements of each sequence. 
\end{lemma}

\begin{defin}[Definition 2.25,~\cite{jikalasa}]
Let $S$ be an unbounded sector of the Riemann surface of the logarithm with vertex at the origin, and let $\mathbb{M}$ be a weight sequence. A function $f$ belongs to the class $\mathcal{O}^{\mathbb{M}}(S)$ if $f$ is continuous at the origin, and for all proper (unbounded) subsector $T$ of $S$ there exist $c,k>0$ such that 
$$|f(z)|\le ce^{\omega_{\mathbb{M}}(|z|/k)},\qquad z\in T.$$ 
\end{defin}

The previous definition can be stated in a more accurate manner, distinguishing different sub-levels in terms of the constant $k>0$.

\begin{defin}
Let $S$ be an unbounded sector of the Riemann surface of he logarithm with vertex at the origin, let $\mathbb{M}$ be a weight sequence, and $k>0$. A function $f$ belongs to the class $\mathcal{O}^{\mathbb{M},k}(S)$ if $f$ is continuous at the origin, and for all proper (unbounded) subsector $T$ of $S$ there exist $c>0$ such that 
$$|f(z)|\le ce^{\omega_{\mathbb{M}}(|z|/k)},\qquad z\in T.$$ 
\end{defin}

\begin{prop}\label{prop00}
Let us suppose that  $f\in\mathcal{O}^{m_{p,q},k}(S)$ for some unbounded sector $S\subset\mathcal{R}$ and some $k>1$. Let $z\in \mathcal{R}$ with $|z|<p/q-1$. Then, there exist a neighborhood $U$ of $z$ and real numbers $\theta,\phi$ which depend on $z$ such that
\begin{equation}\label{e532}
I:=\left|\int_0^{\infty(\theta+ \phi)}\int_0^{\infty(\theta)}\tilde{e}_1(v/z)\tilde{e}_2(u/v)f(u)\frac{du}{u}\frac{dv}{v}\right|
\end{equation}
is finite for all $z\in U$.
\end{prop}
\begin{proof}
Let $r>0$ small enough such that one can choose $\theta,\phi\in\R$ in a way that $\theta\in\hbox{arg}(S)$, $\phi\neq \pi$ and $\phi+\theta+\pi$ does not belong to the arguments of $U=D(z,r)$. From the previous choice, there exist small $\epsilon,\delta>0$ such that if $v\neq0$ with $\hbox{arg}(v)=\phi+\theta$, and $u$ is such that $\hbox{arg}(u)=\theta$, then $|\hbox{arg}(v/z)+\pi|>\epsilon$ together with $|\hbox{arg}(u/v)+\pi|>\delta$. For every fixed $v$ with $\hbox{arg}(v)=\theta+\phi$, we write $v=uw$, with $\hbox{arg}(w)=\phi$. 

Let $s>0$. First, we observe from the definition of the Jacobi Theta function that
$$\Theta_q(1/z)=\Theta_q(z/q).$$
This fact and the estimates in Lemma~\ref{lema1} guarantee that for all $\delta>0$ there exists $\Delta>0$ such that
$$|\Theta_q(qz)|=\left|\Theta_q\left(\frac{1}{q^2z}\right)\right|\ge \Delta \delta \exp\left(\frac{1}{2\log(q)}\log^2\left(\frac{1}{q^2|z|}\right)\right)\frac{1}{q}\frac{1}{|z|^{1/2}},$$
for all $z\in\C\setminus\left(\cup_{m\in\Z}\{z\in\C^{\star}:|1+q^m/z|\le\delta\}\right)$. 

The application of Lemma~\ref{lema4} with $\mathbb{M}_1=m_{p/q}$ and $\mathbb{M}_2=m_1$, yields the existence of $c>0$ such that
\begin{multline*}
\int_0^{s}\left|\tilde{e}_2\left(\frac{r}{se^{i\phi}}\right)\right||f(re^{i\theta})|\frac{dr}{r}=\frac{1}{s}\int_0^{s}\left|e_2\left(\frac{r}{se^{i\phi}}\right)\right||f(re^{i\theta})|dr\\
\le 
\frac{q^{2}}{\log(q)}\frac{1}{\Delta \delta}c e^{\omega_{m_{p/q}}(s)}\frac{1}{s^{3/2}}\int_{0}^{s}\exp\left(-\frac{1}{2\log(q)}\log^2\left(\frac{s}{q^2r}\right)\right)e^{\omega_{m_1}(\frac{r}{ks})}r^{1/2}
dr.
\end{multline*}
The change of variables $s/r=t$ in the integral involving the previous expression leads to an upper bound of the form:
\begin{multline*}
\frac{q^{2}}{\log(q)}\frac{1}{\Delta \delta}c e^{\omega_{m_{p/q}}(s)}e^{\omega_{m_1}(\frac{1}{k})}\int_{1}^{\infty}\exp\left(-\frac{1}{2\log(q)}\log^2\left(\frac{t}{q^2}\right)\right)\frac{1}{t^{5/2}}dt\\
\le \frac{q^{2}}{\log(q)}\frac{1}{\Delta \delta}c e^{\omega_{m_{p/q}}(s)}e^{\omega_{m_1}\left(\frac{1}{k}\right)}\int_{1}^{\infty}\frac{1}{t^{5/2}}dt.
\end{multline*}

Therefore, there exists $\tilde{C}_1>0$ such that
\begin{equation}\label{e567}
\int_0^{s}\left|\tilde{e}_2\left(\frac{r}{se^{i\phi}}\right)\right||f(re^{i\theta})|\frac{dr}{r}\le 
\tilde{C}_1 e^{\omega_{m_{p/q}}(s)}.
\end{equation}

Let us now provide upper estimates for
$$\int_{s}^{\infty}\left|\tilde{e}_2\left(\frac{r}{se^{i\phi}}\right)\right||f(re^{i\theta})|\frac{dr}{r}.$$ 

Taking into account Lemma~\ref{lema1} and Lemma~\ref{lema4} with $\mathbb{M}_1=m_{p/q}$ and $\mathbb{M}_2=m_1$, one arrives at the existence of $\Delta=\Delta(q)>0$ with
$$\int_{s}^{\infty}\left|\tilde{e}_2\left(\frac{r}{se^{i\phi}}\right)\right||f(re^{i\theta})|\frac{dr}{r}\le 
\frac{q^{1/2}}{\log(q)}\frac{1}{\Delta \delta}\frac{1}{s^{1/2}}ce^{\omega_{m_{p/q}}(s)}\int_{s}^{\infty}\exp\left(-\frac{1}{2\log(q)}\log^2\left(\frac{qr}{s}\right)\right)e^{\omega_{m_1}(\frac{r}{ks})}\frac{dr}{r^{1/2}}.$$

Usual computations allow to check (see also 2.6 (4),~\cite{meise}) that 
\begin{equation}\label{e552}
\exp\left(\omega_{m_1}\left(\frac{r}{ks}\right)\right)\le\exp\left(\frac{1}{2\log(q)}\log^2\left(\frac{r}{ks}\right)\right)\left(\frac{r}{ks}\right)^{1/2}C(q),
\end{equation}
for some $C(q)>0$.
This entails that 
\begin{multline*}
\int_{s}^{\infty}\left|\tilde{e}_2\left(\frac{r}{se^{i\phi}}\right)\right||f(re^{i\theta})|\frac{dr}{r}\\
\le \frac{1}{s}\frac{q^{1/2}}{\log(q)}\frac{C(q)}{k^{1/2}\Delta \delta}ce^{\omega_{m_{p/q}}(s)}\int_{s}^{\infty}\exp\left(\frac{1}{2\log(q)}\left(\log^2\left(\frac{r}{ks}\right)-\log^2\left(\frac{qr}{s}\right)\right)\right)dr.
\end{multline*}

From usual estimates one has
\begin{multline*}
\exp\left(\frac{1}{2\log(q)}\left(\log^2\left(\frac{r}{ks}\right)-\log^2\left(\frac{qr}{s}\right)\right)\right)=\exp\left(\frac{1}{2\log(q)}\log\left(\frac{qr^2}{ks^2}\right)\log\left(\frac{1}{kq}\right)\right)\\
=\left(\frac{k}{q}\right)^{\frac{\log(kq)}{2\log(q)}}\left(\frac{s}{r}\right)^{\frac{\log(kq)}{\log(q)}}.
\end{multline*}

This, together with $\int_{s}^{\infty}\frac{dr}{r^{\frac{\log(k)}{\log(q)}+1}}=\frac{\log(q)}{\log(k)}\frac{1}{s^{\frac{\log(k)}{\log(q)}}}$ yields
\begin{equation}\label{e598}
\int_{s}^{\infty}\left|\tilde{e}_2\left(\frac{r}{se^{i\phi}}\right)\right||f(re^{i\theta})|\frac{dr}{r}\le \tilde{C}_2e^{\omega_{m_{p/q}}(s)},
\end{equation}
for some $\tilde{C}_2>0$.

In view of (\ref{e567}) and (\ref{e598}) one arrives at 
$$\int_{0}^{\infty}\left|\tilde{e}_2\left(\frac{r}{se^{i\phi}}\right)\right||f(re^{i\theta})|\frac{dr}{r}\le (\tilde{C}_1+\tilde{C}_2)e^{\omega_{m_{p/q}}(s)},$$
for every $s>0$.

In order to provide an upper bound for $I$ in (\ref{e532}), it is enough to estimate

\begin{equation}\label{e572}
\int_0^{\infty}\left|\tilde{e}_1\left(\frac{se^{i(\theta+\phi)}}{z}\right)\right|e^{\omega_{m_{p/q}}(s)}\frac{ds}{s}=\frac{1}{|z|}\int_0^{\infty}\left|e_1\left(\frac{se^{i(\theta+\phi)}}{z}\right)\right|e^{\omega_{m_{p/q}}(s)}ds.
\end{equation}

From the definition of $e_1$ and the second statement of Proposition~\ref{prop1}, one arrives at

$$\left|\tilde{e}_1\left(\frac{se^{i(\theta+\phi)}}{z}\right)\right|\le \frac{p/q-1}{\log(p/q)}\frac{K_0}{\epsilon}\exp\left(-\frac{\log^2(ps/q|z|)}{2\log(p/q)}\right)\left(\frac{p}{q}\right)^{-\frac{\log(p/q-1)}{\log(p/q)}+\frac{1}{2}}\left(\frac{s}{|z|}\right)^{-\frac{\log(p/q-1)}{\log(p/q)}+\frac{3}{2}},$$
for all $z$ with $|z|\le s(\frac{p}{q}-1)(\frac{p}{q})^{1/2}$. 

We also observe that $q^{\frac{n(n-1)}{2}}\le [n]_q^{!}$, which is a direct consequence of $(q^j-1)/(q-1)\ge q^{j-1}$ for all $j\ge1$. Therefore 
$$e^{\omega_{m_{p/q}}(s)}\le e^{\omega_{m_{1,b}}(s)},$$
for every $s\ge0$, where $m_{1,b}$ stands for the sequence $m_1$ with $q$ substituted by $p/q$. This, together with (\ref{e552}) entail that
$$e^{\omega_{m_{p/q}}(s)}\le \exp\left(\frac{1}{2\log(p/q)}\log^2(s)\right)s^{1/2}C(p/q),$$
for $C(p/q)>0$ as in (\ref{e552}).

Let $H>0$. We first estimate 
$$\int_{H|z|}^{\infty}\left|e_1\left(\frac{se^{i(\theta+\phi)}}{z}\right)\right|e^{\omega_{m_{p/q}}(s)}ds.$$
Regarding the statements above, we observe that the previous expression is bounded from above~by
\begin{equation}\label{e588}
\tilde{C}_3
|z|^{\frac{\log(p/q-1)}{\log(p/q)}-\frac{1}{2}}\int_{H|z|}^{\infty}\exp\left(\frac{1}{2\log(p/q)}(\log^2(s)-\log^2(\frac{ps}{q|z|}))\right)s^{1-\frac{\log(p/q-1)}{\log(p/q)}}ds
\end{equation}
with $\tilde{C}_3=C(p/q)\frac{p/q-1}{\log(p/q)}\frac{K_0}{\epsilon}\left(\frac{p}{q}\right)^{-\frac{\log(p/q-1)}{\log(p/q)}+\frac{1}{2}}$.
Analogous estimates allow us to rewrite (\ref{e588}) in the form
$$\tilde{C}_3
|z|^{\frac{\log(p/q-1)}{\log(p/q)}-\frac{1}{2}}\int_{H|z|}^{\infty}
\exp\left(\frac{1}{2\log(p/q)}\log(\frac{ps^2}{q|z|})\log(\frac{q|z|}{p})\right)
s^{1-\frac{\log(p/q-1)}{\log(p/q)}}ds$$

Assume $|z|<c_0p/q$ for some fixed $c_0\in (0,1)$ such that 
\begin{equation}\label{e643}
c_0<\left(1-\frac{q}{p}\right)\frac{q}{p}.
\end{equation}

We split the previous integral into two parts. 
First, we estimate 
$$\int_{H|z|}^{J|z|^{1/2}}
\exp\left(\frac{1}{2\log(p/q)}\log(\frac{ps^2}{q|z|})\log(\frac{q|z|}{p})\right)
s^{1-\frac{\log(p/q-1)}{\log(p/q)}}ds$$
where $J=(q/p)^{1/2}$. The previous integral reads as follows
\begin{multline*}\int_{H|z|}^{J|z|^{1/2}}
\exp\left(\frac{1}{2\log(p/q)}\log(\frac{ps^2}{q|z|})\log(\frac{q|z|}{p})\right)
s^{1-\frac{\log(p/q-1)}{\log(p/q)}}ds\\
\le \int_{H|z|}^{J|z|^{1/2}}
\exp\left(\frac{1}{2\log(p/q)}\log(\frac{s^2}{c_0})\log(\frac{q|z|}{p})\right)
s^{1-\frac{\log(p/q-1)}{\log(p/q)}}ds\\
\le\int_{H|z|}^{J|z|^{1/2}}
\exp\left(\frac{1}{2\log(p/q)}\log(\frac{H^2|z|^2}{c_0})\log(\frac{q|z|}{p})\right)
s^{1-\frac{\log(p/q-1)}{\log(p/q)}}ds\\
=\int_{H|z|}^{J|z|^{1/2}}
\exp\left(\frac{1}{2\log(p/q)}(\log(\frac{H^2}{c_0})+2\log|z|)(\log(\frac{q}{p})+\log|z|)\right)
s^{1-\frac{\log(p/q-1)}{\log(p/q)}}ds\\
\le\tilde{C}_4\exp\left(\frac{1}{\log(p/q)}\log^2|z|\right)|z|^{\gamma_1},
\end{multline*}
for some $\gamma_1\in\R$, and $\tilde{C}_4>0$.

As for the integral
$$\int_{J|z|^{1/2}}^{\infty}
\exp\left(\frac{1}{2\log(p/q)}\log(\frac{ps^2}{q|z|})\log(\frac{q|z|}{p})\right)
s^{1-\frac{\log(p/q-1)}{\log(p/q)}}ds,$$
it can be estimated from above by
$$\int_{J|z|^{1/2}}^{\infty}
\exp\left(\frac{1}{2\log(p/q)}\log(\frac{ps^2}{q|z|})\log(c_0)\right)
s^{1-\frac{\log(p/q-1)}{\log(p/q)}}ds
$$

$$\le \tilde{C}_5\int_{J|z|^{1/2}}^{\infty}
\exp\left(\frac{1}{2\log(p/q)}\log(\frac{s^2}{|z|})\log(c_0)\right)
s^{1-\frac{\log(p/q-1)}{\log(p/q)}}ds,
$$
for some $\tilde{C}_5>0$. The change of variable $s^2/|z|=t$ turns the previous expression into
$$\tilde{C}_6|z|^{\gamma_2}\int_{J^2}^{\infty}
\exp\left(\frac{\log(c_0)}{2\log(p/q)}\log(t)\right)t^{-\frac{\log(p/q-1)}{2\log(p/q)}}dt,
$$
for $\tilde{C}_6>0$ and some $\gamma_2\in\R$. The previous integral is finite provided that (\ref{e643}) holds.



We proceed to provide upper bounds for
$$\int_{0}^{H|z|}\left|\tilde{e}_1\left(\frac{se^{i(\theta+\phi)}}{z}\right)\right|e^{\omega_{m_{p/q}}(s)}\frac{ds}{s}=\frac{1}{|z|}\int_{0}^{H|z|}\left|e_1\left(\frac{se^{i(\theta+\phi)}}{z}\right)\right|e^{\omega_{m_{p/q}}(s)}ds.$$

The previous integral is bounded from above by
$$e^{\omega_{m_{p/q}}(H|z|)}\int_{0}^{H|z|}\left|e_1\left(\frac{se^{i(\theta+\phi)}}{z}\right)\right|ds.$$
From the definition of $e_1$, and the change of variables $ps/(qz)=t$, one can bound the previous expression from above by
$$\tilde{C}_{7}e^{\omega_{m_{p/q}}(H|z|)}\int_0^{pH/q}\frac{1}{|\exp_{p/q}(te^{i(\theta+\phi-\hbox{arg}(z))})|}dt,$$
for some $\tilde{C}_7>0$. The last integral is convergent in view of the second statement of Proposition~\ref{prop1}.


\end{proof}

As a consequence of the previous result, one arrives at the following integral representation of a kernel function associated with the sequence $m_{p,q}$, together with the associated action of a Laplace-like operator.

\begin{theo}\label{th:3}
The convolution kernel $\tilde{e}$ defined on $\mathbb{R}$ by
$$\tilde{e}(z)=\int_{L_{\phi}}\tilde{e}_1(wz)\tilde{e}_2(1/w)\frac{dw}{w}=z\int_{L_{\phi}}e_1(wz)e_2(1/w)\frac{dw}{w},$$
with $\phi\in\R$, $\phi\neq \pi+2\ell\pi$ for $\ell\in\mathbb{Z}$ satisfies that $\tilde{e}\in\mathcal{O}(\mathcal{R})$ and for all $f\in\mathcal{O}^{m_{p,q},k}(S)$, for some unbounded sector $S\in\mathcal{R}$, and $k>1$, one has that the function $T(f)$ given by
$$T(f)(z)=\int_{L_{\theta}}f(u) \tilde{e}(u/z)\frac{du}{u}$$
defines a holomorphic function for $z\in\mathcal{R}$ with sufficiently small $|z|$. Moreover, 
the $(p,q)-$Gevrey sequence $m_{p,q}$ is the sequence of moments associated with $\tilde{e}$, i.e.,
$$[n]_{p,q}^{!}=\int_0^{\infty}t^{n-1}\tilde{e}(t)dt,\qquad n\ge1.$$
\end{theo}
\begin{proof}
The first part is a direct consequence of Proposition~\ref{prop00}.

The function $T(f)$ is well-defined for $z\in\mathcal{R}$ with $|z|$ small enough as a consequence of Proposition~\ref{prop00}, together with Leibniz's rule and the application of Fubini's Theorem after the change of variable $v=wu$. It is straightforward to check that the construction does not depend on the directions of integration considered from a deformation path argument.

Let us check that the sequence of moments associated with $\tilde{e}$ is the $(p,q)$-Gevrey sequence. Indeed, let $n\ge1$ be an integer. Then, one has that
\begin{multline*}
[n]_{p,q}^{!}=[n]_{p/q}^{!}q^{\frac{n(n-1)}{2}}=\left(\frac{p/q-1}{\log(p/q)}\int_0^{\infty}t^n\frac{1}{\exp_{p/q}(pt/q)}dt\right)\left(\frac{q}{\log(q)}\int_0^{\infty}t^n\frac{1}{\Theta_q(qt)}dt\right)\\
=\frac{p/q-1}{\log(p/q)}\frac{q}{\log(q)}\int_0^{\infty}\int_0^{\infty}(st)^n\frac{1}{\exp_{p/q}(pt/q)}\frac{1}{\Theta_q(qs)}dtds.
\end{multline*}
The change of variables $ts=x_1, 1/s=x_2$ guarantees that the previous expression equals
$$\int_0^{\infty}x_1^n\left(\int_0^{\infty}\frac{p/q-1}{\log(p/q)}\frac{q}{\log(q)}\frac{1}{\exp_{p/q}((p/q)x_1x_2)}\frac{q}{\log(q)}\frac{1}{\Theta_q(q/x_2)}\frac{dx_2}{x_2}\right)dx_1=\int_0^{\infty}x_1^{n-1}\tilde{e}(x_1)dx_1.$$
\end{proof}

\begin{remark} The asymptotic properties on $\tilde{e}$ can be attained by providing more accurate bounds in the previous result, which lead to summability-like results. This is left to a future research which will be devoted to that issue.
\end{remark}

\end{document}